\documentclass{article}
\usepackage{amsmath,amsfonts,amssymb}
\usepackage{authblk}

\usepackage[document]{ragged2e}

\begin{document}
\hspace{2cm} \textbf{An Alternative proof of Steinhaus Theorem} \newline
 \textit{Arpan Sadhukhan, Tata Institute of Fudamental Research, Bangalore, India} \\~\\

Steinhaus's Theorem states that if $A$ is a lebesgue measurable set on the real line such that the lebesgue measure of $A$ is not zero then the difference set $A-A=\{a-b| a,b \in A\}$ contains an open neighbourhood of origin[1]. Since there is a compact set of positive measure inside any set of positive measure. It is enough to prove the result for compact sets. \\~\\

\textbf{Theorem} For any compact set $C\subset \mathbb{R}$ of positive measure, the difference set $D=C-C=\{a-b| a,b \in C\}$ contains an interval containing $0$. \\~\\

\textit{Proof:} Suppose the set $D$ does not contain an interval around origin, then $\exists$ a sequence $x_n \rightarrow 0$ such that $x_n\notin D$, $\forall n\in \mathbb{N}$. \\

We proceed by induction to create arbitrarily large number of mutually disjoint sets $A_1,A_2,A_3\ldots$such that $A_i$ is of the form $x_{n_{i}}+C$ for all $i\geq 2$, $i\in \mathbb{N}$. \newline 

Let $A_1=x_1+C$, clearly $A_1$ satisfies the above property. Now suppose $A_1,A_2\ldots A_k$ is created in such a manner. We will now create $A_{k+1}$. \\

Suppose for every $n \in \mathbb{N}$,  $\exists$ an $i\leq k$ such that the set $x_n+C$ intersects $A_i$, then there exists a $j\leq k$ such that $x_n+C$ intersects $A_j$ for infinitely values of $n\in \mathbb{N}$, so $\exists$ a subsequence $\{y_m\}$ of $\{x_n\}$ such that $x_{n_{j}}-y_m \in D$ for all $m\in \mathbb{N}$. As  $y_m \rightarrow 0$ and $D$ is compact, $x_{n_{j}}$ belongs to $D$ (a contradiction). So there exists an $N\in \mathbb{N}$ such that the set $x_N+C$ does not intersect $A_i$ for any $i\leq k$. Define $A_{k+1}=x_N+C$. So we have our desired set. \newline

Now for any $n\in \mathbb{N}$, the set $x_n+C$ lies in some bounded set $[-R,R]$ for some $R\in \mathbb{N}$ and for any $p,q\in \mathbb{N}$ measure of $A_p>0$, $A_q>0$ and measure of $A_p$ equals measure of $A_q$. So we have an arbitrary large number collection of mutually disjoint sets of the same positive measure in $[-R,R]$ (a contradiction). Hence the result follows. \newline

\section{References}
$[1]$. T.Tao, \textit{An Introduction to Measure Theory}, $1$st edition, American Mathematical Society. \hfill \break \newline

(Accepted American Mathematical Monthly(16/3/19)

\end{document}